\documentclass[12pt]{article}

{\catcode`\~=11
 
 \gdef\picdir{/home/mason/ds/bec/nopot/}
 \gdef\picdir{}
}

\usepackage{fancyheadings,
latexsym,
amsmath,
amstext,
amssymb,
amsgen,
amsbsy,
amsopn,
amsfonts,
graphics,
pstricks,
graphicx,
theorem}

\newtheorem{thm}{Theorem}
\newtheorem{cor}[thm]{Corollary}

\newtheorem{lem}[thm]{Lemma}
\newtheorem{remark}[thm]{Remark}

\newenvironment{rem}{\begin{remark} \rm : }{\end{remark}}
\newcommand{\eq}[1]{(\ref{e:#1})}
\newcommand{\prfend}{\begin{flushright}$\Box$\end{flushright}}
\newenvironment{prf}{{\noindent \bf Proof: }}{\prfend}

\newcommand{\D}{{\cal D}}
\newcommand{\E}{{\cal E}}

\newcommand{\eps}{\varepsilon}
\newcommand{\h}{\frac{1}{2}}

\newcommand{\lnorm}{\left|\left|}
\newcommand{\rnorm}{\right|\right|}

\newcommand{\mC}{{\mathbb C}}

\newcommand{\mQ}{{\mathbb Q}}
\newcommand{\mR}{{\mathbb R}}

\newcommand{\mZ}{{\mathbb Z}}

\newcommand{\re}{\mbox{Re}}
\newcommand{\im}{\mbox{Im}}

\newcommand{\fd}[2]{\frac{\partial #1}{\partial #2}}

\newrgbcolor{lightgray}{.8 .8 .8}

\newcommand{\noteb}[1]{\noindent\textbf{\textit{Authors' note: #1}}}

\newcommand{\change}[1]{#1}
\newcommand{\note}[1]{}

\begin{document}
\title{Quasiperiodic Dynamics in Bose-Einstein Condensates in Periodic Lattices and Superlattices
  }

\author{Martijn van Noort \\ m.vannoort@imperial.ac.uk \\ Department of 
Mathematics \\ Imperial College, London SW7 2AZ, UK\\ \\
  Mason A. Porter\footnote{Corresponding author} \\ mason@caltech.edu \\ 
Department of Physics and Center for the Physics of Information \\
California Institute of Technology, Pasadena, CA  91125 \\ \\
 Yingfei Yi \\ yi@math.gatech.edu \\ Center for Dynamical Systems and Nonlinear
 Studies, School of Mathematics \\ Georgia Institute of Technology, Atlanta, GA  30332 \\ \\ Shui-Nee Chow \\ chow@math.gatech.edu \\ 
Center for Dynamical Systems and Nonlinear Studies, School of Mathematics \\ 
 Georgia Institute of Technology, Atlanta, GA  30332}

\date{\today}

\maketitle

\newpage

\begin{abstract}
  \noteb{We rephrased the abstract.\\}
  We employ KAM theory to rigorously investigate quasiperiodic dynamics in
  cigar-shaped Bose-Einstein condensates (BEC) in periodic lattices and
  superlattices.  Toward this end, we apply a coherent structure
  ansatz to the Gross-Pitaevskii equation to obtain a parametrically
  forced Duffing equation describing the spatial dynamics of the
  condensate.  For shallow-well, intermediate-well, and deep-well
  potentials, we find KAM tori and Aubry-Mather sets to prove that one
  obtains mostly quasiperiodic dynamics for condensate wave functions
  of sufficiently large amplitude, where the minimal amplitude depends
  on the experimentally adjustable BEC parameters.  We show that this
  threshold scales with the square root of the inverse of the two-body scattering length,
  whereas the rotation number of tori above this threshold is
  proportional to the amplitude.  As a consequence, one obtains the
  same dynamical picture for lattices of all depths, as an increase in
  depth essentially only affects scaling in phase space.  Our
  approach is applicable to periodic superlattices with an arbitrary
  number of rationally dependent wave numbers.
\end{abstract}

\vspace{2mm}


\vspace{2mm}


\subsection*{Keywords: Hamiltonian dynamics, Bose-Einstein condensates, KAM 
theory, Aubry-Mather theory}

\section{Introduction}

Bose-Einstein condensates (BECs) have generated considerable
excitement in the physics community both because their study allows
one to explore new regimes of fundamental physics and because of their
eventual engineering applications.  They constitute a macroscopic
quantum phenomenon, and their analysis has already lead to an
increased understanding of phenomena such as superfluidity and
superconductivity.  Of particular interest are BECs in optical
lattices (periodic potentials), which have already been used to study
Josephson effects \cite{anderson}, squeezed states \cite{squeeze},
Landau-Zener tunneling and Bloch oscillations \cite{morsch}, and the
transition between superfluidity and Mott insulation \cite{smerzi,cata,mott}. With each lattice site occupied by one alkali atom in its ground state, BECs in periodic potentials also
show promise as registers for quantum computers \cite{porto,voll}.

\noteb{We rephrased the next paragraph.}

In the present paper, we generalize recent work on near-autonomous
dynamics in BECs \cite{mapbecprl,mapbec} to study quasiperiodic behavior in
BECs in periodic lattices, which can have shallow, intermediate, or
deep wells.  We present our methodology and results in section
\ref{results}.  In section \ref{setup}, we discuss the physics of BECs
and use a coherent structure ansatz to derive a parametrically forced
Duffing oscillator describing the spatial dynamics of the condensate.
The quasiperiodic dynamics of parametrically forced Duffing oscillators is
rigorously investigated in section \ref{main}.  Sections \ref{lattice} and
\ref{superlattice} then describe applications to BECs in periodic
lattices and periodic superlattices, respectively. The KAM theorem
used in this analysis is proven in section \ref{s:proof}.  Finally, we summarize and discuss our results in section \ref{conc}.

\subsection{Methodology and results}\label{results}

\noteb{We rephrased the paragraph below.}

The spatial dynamics of standing waves in BECs in periodic
optical lattices can be described by a parametrically forced Duffing
equation, where the periodic forcing is given by an external potential
due to the lattice \cite{bronski,mapbecprl,mapbec}.  This gives a $1\h$ degree of freedom Hamiltonian system.  We use KAM theory and Aubry-Mather theory to study its quasiperiodic dynamics.
 Previous KAM studies in BECs took a heuristic approach and considered only near-autonomous situations \cite{mapbec,mapbecprl}.  The approach of this paper, however, is especially versatile in that shallow, deep, and intermediate lattice wells can all be considered using the same mathematical framework.  That is, we consider the near-autonomous and far-from-autonomous settings simultaneously.

Theorem \ref{t:main} proves that for any (analytic) external periodic
potential, any negative two-body scattering length, and any chemical potential,
one obtains mostly 2-quasiperiodic dynamics for condensate wave
functions of sufficiently large amplitude, where the minimal amplitude
depends on the experimentally adjustable BEC parameters.  In
particular, the threshold amplitude is proportional to the reciprocal
of the square root of the scattering length. Any 2-quasiperiodic wave
function above the threshold has one fixed frequency and one
proportional to its amplitude.  We also demonstrate numerically that
one obtains the same dynamical picture for lattices of all depths, as
an increase in lattice amplitude essentially only affects scaling in
phase space.  These numerical results support the theoretically
predicted scaling of the threshold amplitude. Our theorem applies to
periodic superlattices with an arbitrary number of rationally
dependent wave numbers.

The system we investigate is given by a Hamiltonian of the form
$H(R,S,\xi) = \h S^2 + U(R,\xi)$, where $|R|$ is the amplitude of the
wave function, $R$ and $S$ are conjugate variables, $\xi \in \mR/\mZ$, and $\xi'=1$.  The function $U$ is a polynomial in $R$ and 1-periodic in $\xi$.
We also consider its Poincar\'e map, defined to be the return map on
the section $\xi = 0$.

One can show that such systems have invariant tori sufficiently far from the
origin $(R,S) = (0,0)$ even when they are far from autonomous (that is, even
when $H$ is a large perturbation of an $\xi$-independent system). In the
present case, this implies the existence of invariant tori for any optical
lattice depth. The key condition is that $U(R,\xi)/R^2 \to +\infty$ as $R \to
\pm \infty$, as this guarantees that the set of frequencies corresponding to
rotation around the origin is unbounded; indeed, the frequency goes to
infinity with the distance to the origin. The first result of this type can be
found in the equivalent context of adiabatic theory \cite{Ar63b}.  Another
strand goes back to the question of boundedness of solutions (which is implied
by the existence of invariant tori in this low-dimensional situation)
\cite{Litt66}. These qualitative results have also been extended to more
general mathematical settings \cite{DZ87,Mo89a,Mo89b,Liu89,You90,Le91,LZ95}.

\noteb{We added Levi91 to the above list of references.}

\begin{figure}
  \centerline{(a) \hspace{0.285 \textwidth} (b) \hspace{0.285 \textwidth} (c)}
  \centerline{
    \includegraphics[width=0.3\textwidth, height=0.3\textwidth]
    {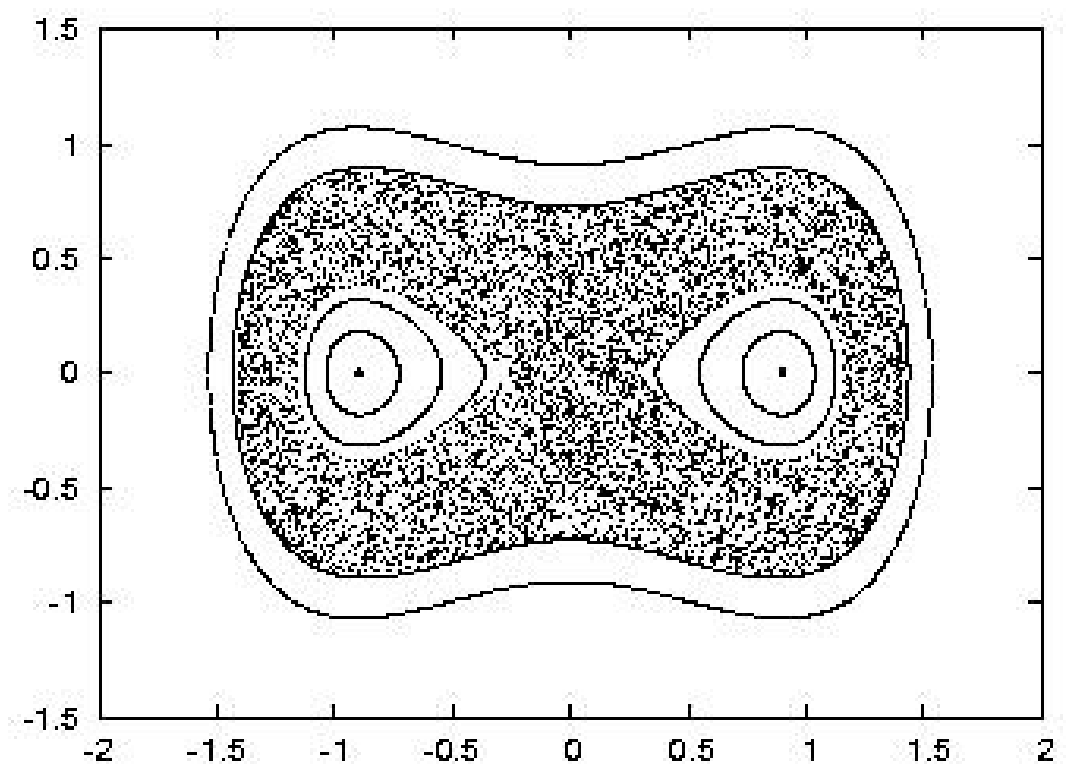}
    \hspace{0.02 \textwidth}
    \includegraphics[width=0.3\textwidth, height=0.3\textwidth]
    {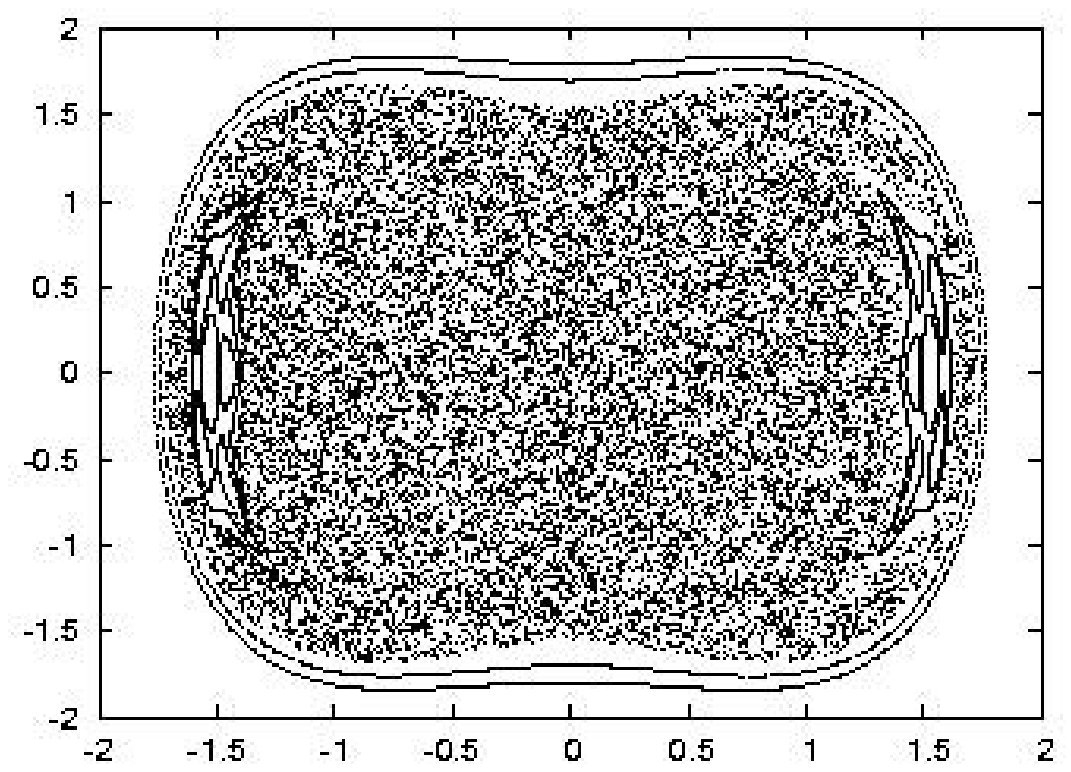}
    \hspace{0.02 \textwidth}
    \includegraphics[width=0.3\textwidth, height=0.3\textwidth]
    {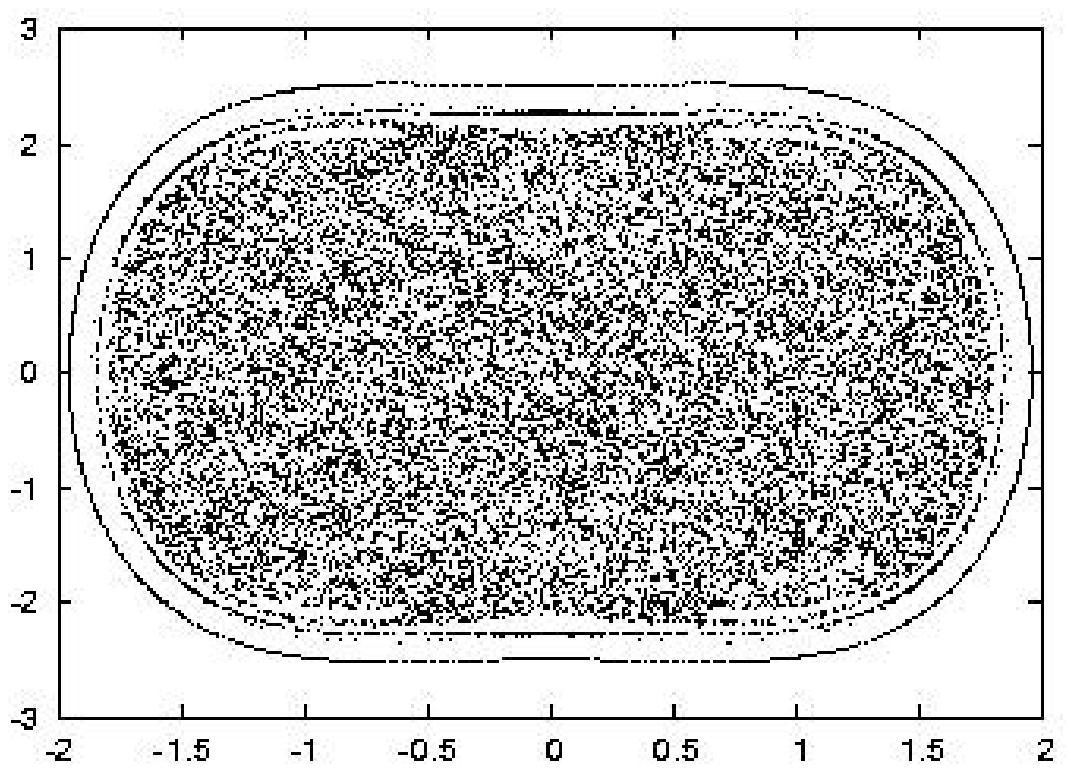}}
  \caption{Phase portraits of the Poincar\'e map $P$ for the example
    of section \ref{lattice}, with horizontal axis $R$ and vertical axis $S$. In this case, the Hamiltonian is of the form $H = \h S^2 - \h [1+V_1 \cos(\xi)] R^2 + \frac{1}{4}
    R^4$.  Observe the invariant curves for large $R$
    in these figures.  Generically, the invariant manifolds of the
    central saddle point intersect transversally, creating a
    homoclinic tangle and thereby implying the existence of horseshoes of
    measure zero. By conjecture, the closure of these horseshoes is a
    ``chaotic sea'' of positive measure, corresponding to the dots in
    the figures.  As the size of the perturbation $V_1$ is increased
    to $+\infty$, the size of each of the remaining integrable islands
    vanishes, but their total measure remains $O(1)$. The values of
    $V_1$ are (a) $V_1 = 0.1$, (b) $V_1 = 0.5$, and (c) $V_1 = 1$.}
  \label{p:pp1}
\end{figure}

\begin{figure}
  \centerline{(d) \hspace{0.285 \textwidth} (e) \hspace{0.285 \textwidth} (f)}
  \centerline{
    \includegraphics[width=0.3\textwidth]
    {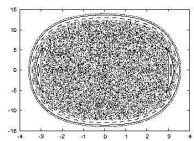}
    \hspace{0.02 \textwidth}
    \includegraphics[width=0.3\textwidth]
    {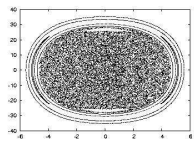}
    \hspace{0.02 \textwidth}
    \includegraphics[width=0.3\textwidth]
{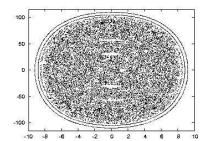}}
  \caption{Continuation of figure \ref{p:pp1}. Phase portraits of the
    Poincar\'e map $P$ at (d) $V_1 = 10$, (e) $V_1 = 25$, and (f) $V_1
    = 100$.}
  \label{p:pp2}
\end{figure}


\begin{figure}
  \centerline{(g) \hspace{0.285 \textwidth} (h) \hspace{0.285 \textwidth} (i)}
  \centerline{
    \includegraphics[width=0.3\textwidth]
    {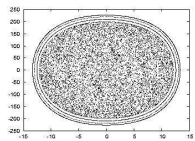}
    \hspace{0.02 \textwidth}
    \includegraphics[width=0.3\textwidth]
    {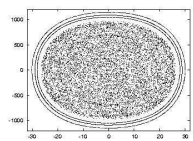}
    \hspace{0.02 \textwidth}
    \includegraphics[width=0.3\textwidth]
    {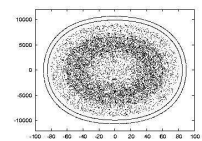}}
  \caption{Continuation of figure \ref{p:pp2}. Phase portraits of the
    Poincar\'e map $P$ at (g) $V_1 = 200$, (h) $V_1 = 1000$, and (i) $V_1
    = 10000$.}
  \label{p:pp3}
\end{figure}

\noteb{The paragraph below was rephrased.}

A typical system of this type exhibits a phase space divided into two
clearly distinct regions.  See figures \ref{p:pp1} -- \ref{p:pp3},
which show numerical experiments in the present setting. One region,
bounded away from the origin, consists largely of invariant tori, with
small layers of chaotic dynamics between them. Indeed, one can show
that the measure of these layers vanishes exponentially fast as the
distance to the origin goes to infinity \cite{Nei81,BNS03}.  The other region has mainly chaotic dynamics.  In this paper, we employ a quantitative existence result \cite{cny}
to obtain bounds on the location and frequencies of invariant tori.

\section{Physical Background}\label{setup}

At low temperatures, particles in a dilute boson gas can reside in the
same quantum (ground) state, forming a Bose-Einstein condensate (BEC)
\cite{stringari}.  This was first observed experimentally in 1995 with
vapors of rubidium and sodium \cite{becrub,becna}. In these
experiments, atoms were confined in magnetic traps, evaporatively
cooled to temperatures of a few hundred nanokelvin, left
to expand by switching off the confining trap, and subsequently imaged
with optical methods \cite{stringari}. A sharp peak in the velocity
distribution was observed below a critical temperature, indicating
that condensation had occurred.  BECs are inhomogeneous, allowing
condensation to be observed in both momentum and position space.  The
number of condensed atoms $N$ ranges from several thousand to tens of millions.

A BEC has two characteristic length scales: the harmonic oscillator
length $a_{ho} = \sqrt{\hbar/[m\omega_{ho}]}$ (which is about a few
microns), where $\omega_{ho}=(\omega_x \omega_y \omega_z)^{1/3}$ is
the geometric mean of the trapping frequencies, and the mean healing
length $\chi=1/\sqrt{8\pi |a| \bar{n}}$ (which is also about a few microns), where $\bar{n}$ is the mean
density and $a$, the (two-body) $s$-wave scattering length, is
determined by the atomic species of the condensate.  Interactions
between atoms are repulsive when $a > 0$ and attractive when $a < 0$.
For a dilute ideal gas, $a \approx 0$. The length scales in BECs
should be contrasted with those in systems like superfluid helium, in
which the effects of inhomogeneity occur on a microscopic scale fixed
by the interatomic distance \cite{stringari}.

\note{We rephrased the paragraph below slightly to discuss the physics a bit more.}

\change{When considering only two-body interactions, the BEC wave function 
(``order parameter'') $\Psi(\vec{r},t)$ satisfies the Gross-Pitaevskii (GP) equation,
\begin{equation}
  	i\hbar {\Psi }_{t}=\left( -\frac{\hbar ^{2}\nabla ^{2}}{2m}+g_{0}|\Psi|^{2}+{\mathcal{V}}(\vec{r})\right) \Psi\,,  \label{GPE}
\end{equation}
where $\Psi =\Psi (\vec{r},t)$ is the condensate wave function, ${\mathcal{V}}(\vec{r})$ is the external potential, and the effective interaction constant is $g_{0}=[4\pi\hbar ^{2}a/m][1+O(\zeta ^{2})]$, where $\zeta \equiv \sqrt{|\Psi |^{2}|a|^{3}}$ is the dilute-gas parameter \cite{stringari,lieb}.  A BEC is modeled in the quasi-one-dimensional (quasi-1D) regime when its transverse dimensions are
on the order of its healing length and its longitudinal dimension is
much larger than its transverse ones \cite{bronski,stringari}.  In the quasi-1D
regime, one employs the 1D limit of a 3D mean-field theory rather than
a true 1D mean-field theory, which would be appropriate were the
transverse dimension on the order of the atomic interaction length or
the atomic size.  The resulting 1D equation is \cite{stringari,salasnich}
\begin{equation}
        i\hbar\psi_t = -\frac{\hbar^2}{2m}\psi_{xx} + g|\psi|^2\psi + V(x)\psi 
\,, \label{nls3}
\end{equation}
where $\psi$, $g$, and $V$ are, respectively, the rescaled 1D wave function, interaction constant, and external trapping potential.   The quantity $|\psi|^2$ gives the atomic number density.  The self-interaction parameter $g$ is tunable (even its sign), because the scattering length $a$ can be adjusted using magnetic fields in the vicinity of a Feshbach resonance \cite{fesh,FRM}.} 

Potentials $V(x)$ of interest include harmonic traps, periodic
lattices and superlattices (i.e., optical lattices with two or more
wave numbers), and periodically perturbed harmonic traps.  The
existence of quasi-1D cylindrical (\textquotedblleft cigar-shaped'') BECs motivates
the study of periodic potentials without a confining trap along the
dimension of the periodic lattice \cite{band}.  Experimentalists use a
weak harmonic trap on top of the periodic lattice or superlattice to
prevent the particles from spilling out.  To achieve condensation, the
lattice is typically turned on after the trap.  If one wishes to
include the trap in theoretical analyses, $V(x)$ is modeled by
\begin{equation}
        V(x) = V_1\cos(\kappa_1 x) + V_2\cos(\kappa_2 x) + V_h x^2 \,, 
\label{harmwiggle}
\end{equation}
where $\kappa_1$ is the primary lattice wave number, $\kappa_2 >
\kappa_1$ is the secondary lattice wave number, $V_1$ and $V_2$ are
the associated lattice amplitudes, and $V_h$ represents the magnitude
of the harmonic trap.  (Note that $V_1$, $V_2$, $V_h$, $\kappa_1$, and
$\kappa_2$ can all be tuned experimentally.)  When $V_h \ll
V_1\,,V_2$, the potential is dominated by its periodic contributions
for many periods.  BECs in optical lattices with up to $200$ wells
have been created experimentally \cite{well}.

In this work, we let $V_h = 0$ and focus on periodic lattices and
superlattices.  Spatially periodic potentials have been employed in
numerous experimental studies of BECs (see, for example, Refs~\cite{anderson,hagley}) and have also been studied theoretically (see, for example, \cite{bronski,promislow,pethick2,wu4}).  In recent experiments, BECs have also been loaded successfully into superlattices with $\kappa_2
= 3\kappa_1$ \cite{quasibec}.  Additionally, over the past couple years, there has been an increasing number of theoretical studies on BECs in superlattices \cite{eksioglu,louis2,super,superlocal,chua}.

\noteb{The next paragraph has been rephrased.}

As mentioned in section \ref{results} and proven below, we obtain bounds on the location and frequencies of invariant tori.  Our paper establishes quasiperiodic dynamics for sufficiently large amplitude.
By a slight adaptation of the argument (see remark \ref{t:measure} below), it can be shown that the quasiperiodic dynamics has large measure at large amplitude. It is generally conjectured that
chaotic dynamics exist where invariant tori have been destroyed.
In the present model, the presence of chaos (which was studied in BECs in, e.g., Refs.~\cite{promislow,chong,thom,chua}) would reflect an irregular spatial profile $R(x)$ (where $|R(x)| = |\psi(x,t)|$ is the amplitude of the BEC wave function), whereas invariant tori correspond to regular (i.e., quasiperiodic) spatial profiles.

{\bf map: how is the phrasing above?}

\begin{rem}
  When the optical lattice has deep wells (large $|V_1|$ or $|V_2|$),
  one can also obtain an analytical description of BECs in terms of
  Wannier wave functions using the so-called ``tight-binding
  approximation'' \cite{smer2}.  In this regime, the BEC dynamics is
  governed by a discrete nonlinear Schr\"odinger equation, which is
  derived by expanding the field operator in a Wannier basis of
  localized wave functions at each lattice site.
\end{rem}


Coherent structures solutions are described with the ansatz 
\begin{equation}
	  \psi(x,t)= R(x)\exp\left(i\left[\theta(x) - \mu t\right]\right)\,, 
\label{maw2} 
\end{equation}
where $R \in \mathbb{R}$ gives the amplitude dynamics of the wave
function, 
$\theta$ gives the phase dynamics, and the ``chemical potential'' $\mu$, defined as the energy
it takes to add one more particle to the system, is proportional to
the number of atoms trapped in the condensate.
When the (temporally periodic) coherent structure (\ref{maw2}) is also
spatially periodic, it is called a {\it modulated amplitude wave} (MAW) \cite{mapbec}.

\begin{rem}
The present work is concerned with the spatial amplitude dynamics of solutions
of the form (\ref{maw2}).  To ensure that such solutions are physically
relevant, it is important to examine their stability with respect to the
dynamics of the GP equation (\ref{nls3}).  Bronski and coauthors
\cite{bronskiatt,bronski,bronskirep} were able to obtain rigorous stability
results in some situations using elliptic-function potentials.  (For repulsive
condensates, for example, they used the ansatz (\ref{maw2}) to construct elliptic-function solutions to
(\ref{nls3}), whose linear stability they proved for $R(x) > 0$ by showing that they were ground states of the GP equation.)  Their results can be applied to trigonometric potentials by taking the limit as the elliptic modulus approaches zero.  More generally, one can address the stability of solutions of the form (\ref{maw2}) through direct numerical simulations of (\ref{nls3}) using such solutions as initial wave functions: $\psi(x,0) = R(x)$.  It was shown previously using both lattice and superlattice potentials $V(x)$ that one can obtain numerically stable MAWs for (\ref{nls3}) with solutions of the form (\ref{maw2}) with trivial phase $\theta = 0$ \cite{super,mapbin}.  Furthermore, the numerically stable \textquotedblleft period-doubled" solutions (whose spatial periodicity is twice that of the optical lattice potential) constructed via subharmonic resonances using the ansatz (\ref{maw2}) with trivial phase ($\theta = 0$) have very recently been observed experimentally \cite{chu}.
\end{rem}

Inserting (\ref{maw2}) into the GP equation (\ref{nls3}) and equating real and
imaginary parts, one obtains
\begin{align}
        \hbar\mu R(x) &= -\frac{\hbar^2}{2m}R''(x)
        + \left[\frac{\hbar^2}{2m}\left[\theta'(x)\right]^2 + gR^2(x) 
+ V(x) \right]R(x) \,, \\
        0 &= \frac{\hbar^2}{2m}\left[2\theta'(x)R'(x) + \theta''(x)R(x)\right] 
             \,, \notag
\end{align}
which gives the following nonlinear ordinary differential equation:
\begin{equation}
  R'' = \frac{c^2}{R^3} - \frac{2m\mu R}{\hbar} 
+ \frac{2mg}{\hbar^2}R^3 + \frac{2m}{\hbar^2}V(x)R \,. \label{moment}
\end{equation}
The parameter $c$ is defined via the relation
\begin{equation}
        \theta'(x) = \frac{c}{R^2(x)}\,,  \label{angmom}
\end{equation}
which plays the role of conservation of ``angular momentum,'' as 
discussed by Bronski and coauthors \cite{bronski}.  Constant phase solutions 
(i.e., standing waves) constitute an important special case and satisfy 
$c = 0$.  In the rest of the paper, we consider only standing waves, so that
\begin{equation}
  R'' = - \frac{2m\mu R}{\hbar} 
  + \frac{2mg}{\hbar^2}R^3 + \frac{2m}{\hbar^2}V(x)R \,. \label{dynam35}
\end{equation}

\begin{figure}
                \centerline{
                (a)
                \includegraphics[width=0.3\textwidth]
{\picdir 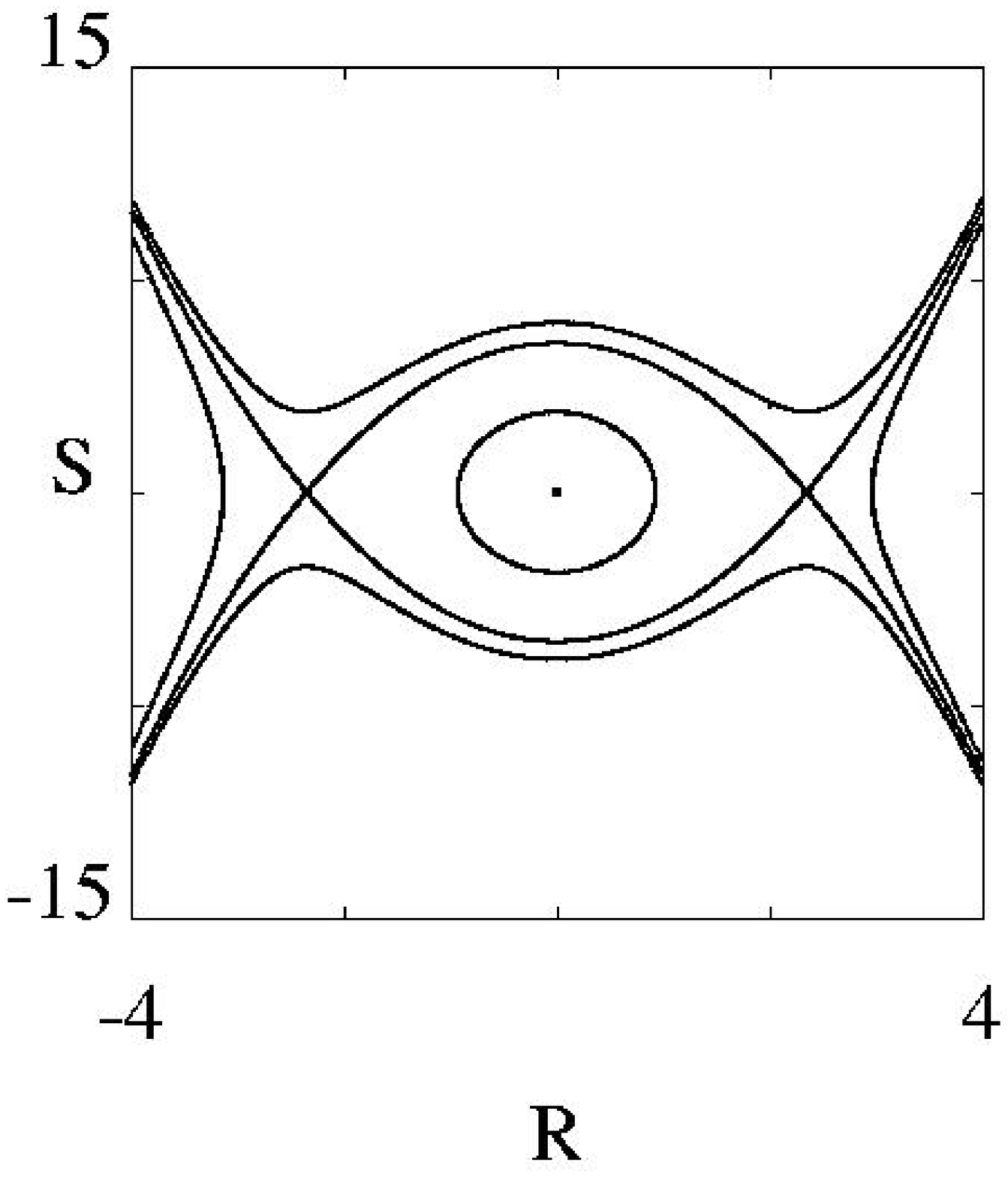}
                (b)
                \includegraphics[width=0.3\textwidth]
{\picdir 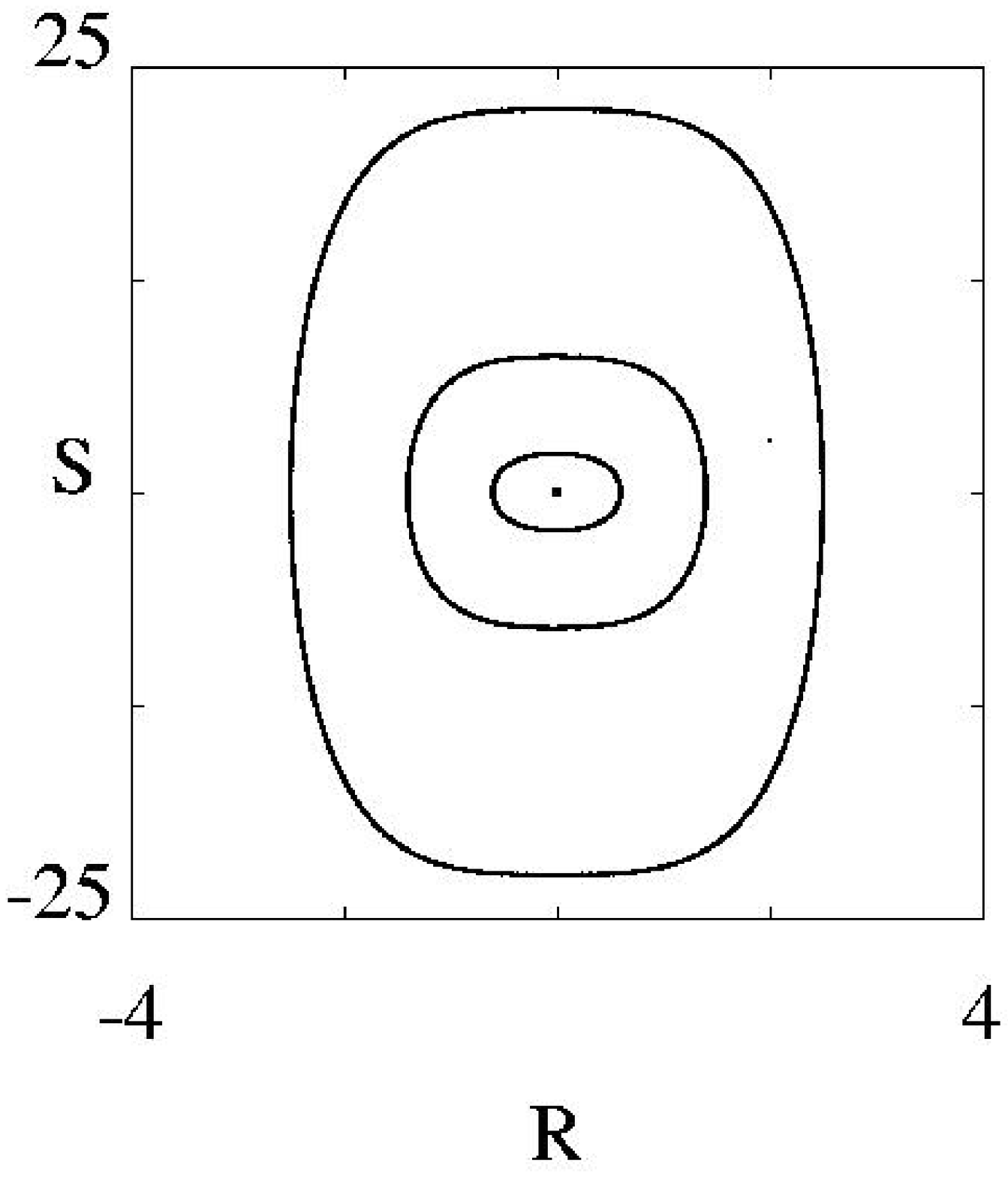}
                (c)
                \includegraphics[width=0.3\textwidth]
{\picdir 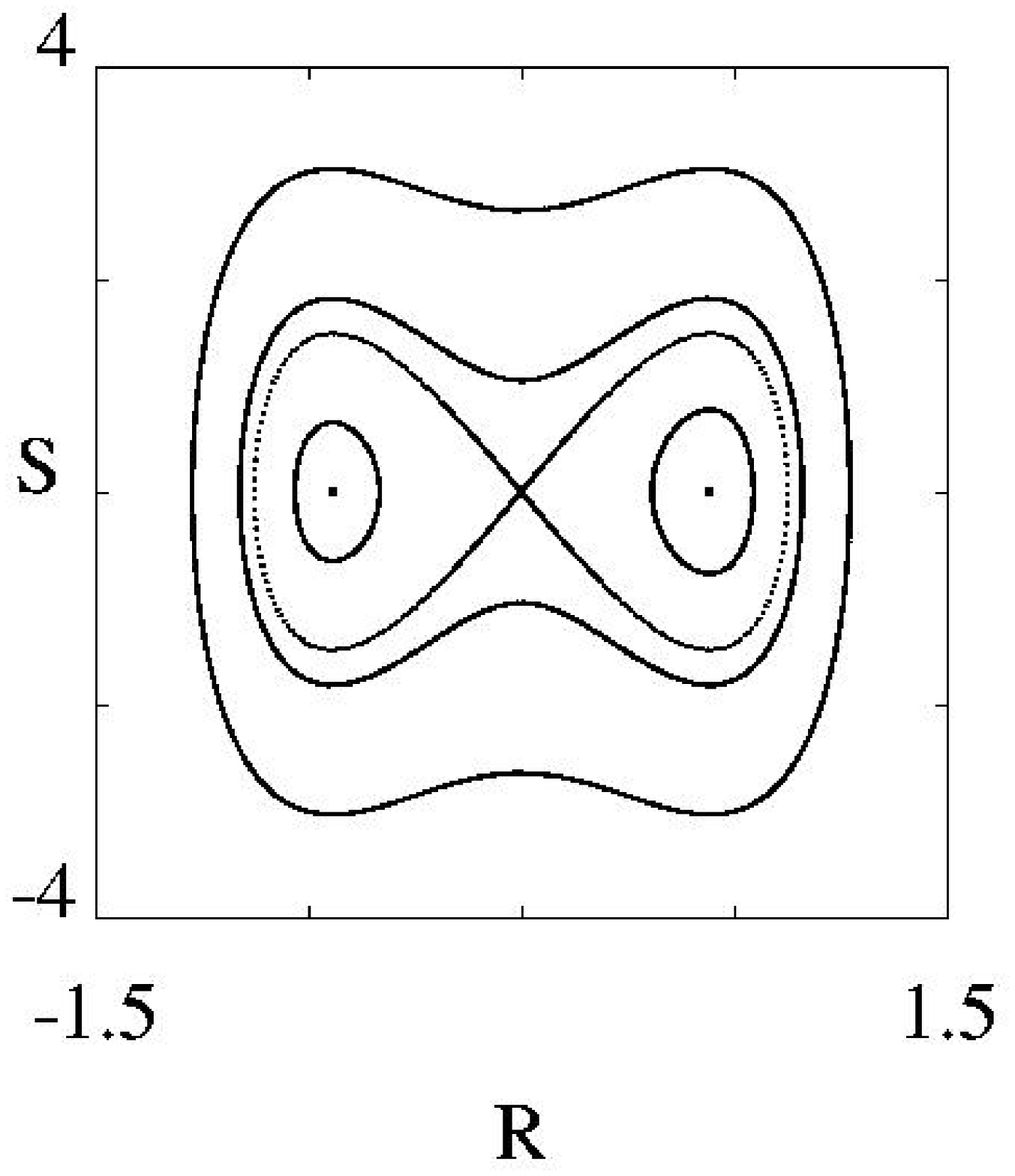}}


                \caption{Phase portraits of coherent structures in BECs with 
                  no external potential.  The signs of $\mu$ and $g$
                  determine the dynamics of (\ref{dynam35}).  (a)
                  Repulsive BEC with $\mu > 0$.  Orbits inside the
                  separatrix have bounded amplitude $|R(x)|$.  The
                  period of such orbits increases as one approaches
                  the separatrix.  In this case, the dynamical system
                  can be rescaled so that $R'' = -R + R^3$.  (b)
                  Attractive BEC with $\mu > 0$.  The dynamical system
                  can be rescaled so that $R'' = -R - R^3$.  (c)
                  Attractive BEC with $\mu < 0$.  Here there are two
                  separatrices, each of which encloses periodic orbits
                  satisfying $R \neq 0$.  The dynamical system can be
                  rescaled so that $R'' = R - R^3$.} \label{repulse1}
\end{figure}

\begin{rem}
  When $V(x) \equiv 0$, the dynamical system (\ref{dynam35}) is the
  autonomous, integrable Duffing oscillator. Its qualitative dynamics
  in the physically relevant situation of bounded $|R|$ is
  illustrated in figure \ref{repulse1}.  The methodology developed in
  the present paper can handle attractive BECs ($g < 0$) with either
  $\mu < 0$ or $\mu > 0$ but not repulsive BECs, as equation
  (\ref{dynam35}) has unbounded solutions when $g > 0$.
\end{rem}


\section{Main Result}\label{main}

This section states the main theorem of this paper, which concerns the
existence of quasiperiodic dynamics in a class of systems of the form
(\ref{dynam35}), including the cases with $V(x)$ given by periodic lattices and
superlattices. Applications to these cases are subsequently given in
two subsections.

The equation of motion (\ref{dynam35}) describes a $1\h$ degree of freedom
Hamiltonian system whose forcing is periodic. We will allow arbitrary analytic
periodic potential functions $V$.  We first rescale the period of the forcing to 1.  We introduce the
phase variable $\xi = T^{-1} x \,(\mbox{mod 1})$, where $T > 0$ is the
minimal period of $V$. Letting $'$ denote differentiation with
respect to $\xi$, we define $S = R'$, and
\[ z_2(\xi) = T^2 \left( \frac{m \mu}{\hbar} - \frac{m}{\hbar^2} V(x)
\right)\,, \quad z_4 = -T^2 \frac{m g}{2 \hbar^2}\,, \quad U(R,\xi) =
z_2(\xi) R^2 + z_4 R^4. \]
This gives the suspended dynamical system
\begin{align}
  R' &=  S\,, \notag \\
  S' &= -\fd{U}{R}(R,\xi)\,, \notag \\
  \xi' &= 1 \,,\label{e:sys}
\end{align}
with Hamiltonian
\begin{equation}
  H(R,S,\xi) = \h S^2 + U(R,\xi)\,. \label{e:H}
\end{equation}

\noteb{We changed the following remark.}

\begin{rem}
  We restrict to analytic systems for convenience only.  It is not a
  necessary restriction, as our result is based on Herman's
  \emph{translated curve theorem}~\cite{Her83,Her86}, which requires only $C^4$
  smoothness. Furthermore, finite-smoothness KAM results can
  be obtained from analytic ones by inverse approximation \cite{Poe82b}.
  The restriction to analytic systems allows one to use Cauchy's integral
  theorem to replace $C^k$ norms by $C^0$ ones.
\end{rem}

\begin{rem}
  Because $z_4 > 0$, the function $U$ is a well; that is, $U \to
  +\infty$ as $R \to \pm \infty$ for all $\xi$.  It is even in $R$,
  which has important ramifications for the sizes of the perturbations
  in our subsequent analysis.  (See, for example, lemma
  \ref{t:sys_rescaled} below, where the leading term in the
  nonintegrable parts $F_1$ and $F_2$ would have been order one
  instead of going to zero as $R \to +\infty$ were it not for this
  symmetry.)
\end{rem}

We introduce action-angle coordinates as follows.  (The details are in
section \ref{s:proof}.)  Let $H_0(R,S) = \h S^2 + z_4 R^4$. For $h >
0$, define the action $I = I(h)$ to be the area in the $(R,S)$-plane
enclosed by the curve $H_0(R,S) = h$. Let the angle $\phi = \phi(R,S)
\in \mR/\mZ$ be such that the transformation $(R,S) \mapsto (\phi,I)$ is
symplectic. This defines $\phi$ uniquely if we set $\phi(0,S) \equiv 0$ for $S > 0$.

In action-angle coordinates, the Hamiltonian takes the form
$K(\phi,I,\xi) = K_0(I) + K_1(\phi,I,\xi)$, where $K_0(I) = H_0(R,S)$.
We consider $K$ as a perturbation of $K_0$. For any $I_0 > 0$, the
unperturbed system $K_0$ has an invariant torus $I=I_0$ with frequency
$\omega = K_0'(I_0)$ in $\phi$.  We say that this frequency is of {\it
  constant type} with parameter $\gamma > 0$ if
\begin{equation}\label{e:dioph}
  \left|\omega - \frac{p}{q}\right| \geq \gamma q^{-2} \quad \mbox{
    for all} \quad \frac{p}{q} \in \mQ\,.
\end{equation}
This is a special type of Diophantine condition \cite{Her83,Her86}.

For a function $f$ defined on a set $\D$, we define $||f||_{\D} =
\sup_{\D} |f|$. If $f$ is vector-valued, then $||f||_{\D}$ is the
maximum of the norms of the components. For $d > 0$, let $\bar{\D}(d) =
\{ \xi \in \mC/\mZ : |\im(\xi)| \leq 2 d \}$.

With the additional notation $\eta = 18 \gamma$, $M = z_4 I_0$,
$c=126/25 = 5.04$, and $b_1 = \int_0^1 (1-u^4)^{1/2} \, {\rm d} u =
0.874019\ldots$, the main result of this paper can now be stated as
follows.
\begin{thm}\label{t:main}
  The Hamiltonian $K$ has an invariant torus with frequency $\omega =
  K_0'(I_0)$ in $\phi$ if there exist $\nu$, $\gamma$, $d > 0$, and
  $b_2 > 0$ with $0 \leq \nu \leq \frac{1}{9} \cdot 2^{-7/3}$ and $0 <
  \gamma \leq \frac{49}{72}$ such that $\omega$ is of constant type
  with parameter $\gamma$ and the following conditions hold:
  \begin{eqnarray}
    A \left(1 + \frac{3 \log{B}}{\log{M}}\right) & \leq &
    2^{-4/3} b_1^{-4/3} \log(2)\,, \label{e:cond1}\\
    1 & < & M\,, \label{e:condM}\\
    L & \leq & \frac{47}{200}\,, \label{e:condL}\\
    b_2 & \leq & 18 [1+24 b_1 (\eta+2 d)] (1+L)^{2/3} (1-L)^{1/3}
    \log(M)\,, \label{e:condF}\\
    2 & \leq & B M^{1/3}\,, \label{e:condb}\\
    ||z_2||_{\bar{\D}(d)} & \leq & \delta \frac{M^{2/3}}{\log^2(M)}\,,
    \label{e:cond_last}
  \end{eqnarray}
  where
  \begin{eqnarray*}
    L & = & 2^{4/3} \cdot 3 b_1^{4/3} \eta M^{-1/3} + \frac{2 b_2
      d}{\log(M)}\,,\\
    A & = & \frac{1 + 4 (\eta+2 d)}{d}
    \max\left\{ \frac{2 \eta}{3} M^{-1/3} \log(M) +
      \left( \frac{2^{2/3}}{9} + \frac{7 \nu}{2}\right) d b_1^{-4/3}
      b_2 +\right.\\
    & & \left.\frac{2^{-1/3}}{27} b_1^{-4/3} (1-L)^{-5/3} L^2 \log(M),
      2 M^{-1/3} \log(M) + \frac{3 \nu}{2} d b_1^{-4/3} b_2 \right\}\,,\\
    B & = & 2^{-13/3} 3^3 7^8 c^{-1} d b_1^{-4/3} b_2 \eta^{-6}
    \left(\frac{1}{\log(M)} + \frac{2^{-7/3}}{3} b_1^{-4/3} b_2
      \frac{M^{1/3}}{\log^2(M)} \right)\,,\\
    \delta & = & \frac{2^{-1/3}}{3} \nu d b_2^2 b_1^{-5/3}
    (1+L)^{-2/3} [1+24 b_1 (\eta+2 d)]^{-3}.
  \end{eqnarray*}
  The torus lies in the region given by $|I-I_0| \leq b_2 I_0
  [\log(z_4 I_0)]^{-1} (\frac{11}{19} d + \rho)$.
\end{thm}

\begin{rem}
  Conditions \eq{cond1} -- \eq{condb} are satisfied for $M$
  sufficiently large, while \eq{cond_last} is a
  restriction on $z_2$ and hence on $V$.
  The theorem implies that $||z_2||_{\bar{D}(d)}$ can be taken roughly
  proportional to $z_4^{2/3} I_0^{2/3}$. Equation \eq{I} below shows
  that $z_4^{-1/6} I_0^{1/3}$ is proportional to the maximal $R$
  coordinate $R_{\max}$ on the torus $H_0(R,S) = K_0(I_0)$, so
  $||z_2||_{\bar{D}(d)}$ is roughly proportional to $z_4 R_{\max}^2$. This
  fits with the numerically computed phase  portraits in figures
  \ref{p:pp1} -- \ref{p:pp3}.
  
  In terms of physical parameters, this implies that for an attractive
  BEC with given scattering length $a < 0$ and chemical potential
  $\mu$ loaded into an arbitrary periodic lattice of amplitude $||V||$
  (with any number of wave numbers), the wave function's spatial
  component $R(x)$ is quasiperiodic with two frequencies 
  if its maximum $R_{\max}$ is large enough, its frequencies
  satisfy the Diophantine condition \eq{dioph}, and the amplitude
  $||V||$ of the lattice potential is sufficiently small.  The lower
  bound on $R_{\max}$ scales as $\kappa/\sqrt{|a|}$, where $\kappa =
  2\pi/T$ is the lattice wave number (and we recall that $a$ is the
  two-body scattering length), whereas the upper bound on
  $||V||$ scales as $a R_{\max}^2$. All frequency ratios that are
  algebraic numbers of index 2 satisfy the Diophantine condition (for
  some $\gamma > 0$).
\end{rem}

\begin{rem}
  The conditions of the theorem imply that $M$ should be larger than
  roughly $10^6$. Indeed, from $A \geq 8 \cdot 2 M^{-1/3} \log(M)$ and
  $B M^{1/3} \geq 2$, it follows that $A [1 +
  3 \log(B)/\log(M)] \geq 48 M^{-1/3} \log(2)$. Hence condition
  \eq{cond1} implies that $M \geq 2^{16} 3^3 b_1^4 \approx 10^6$.
\end{rem}

\begin{figure}[htp]
  \begin{picture}(150,70)
    \put(25,0){\resizebox{10cm}{!}{\includegraphics{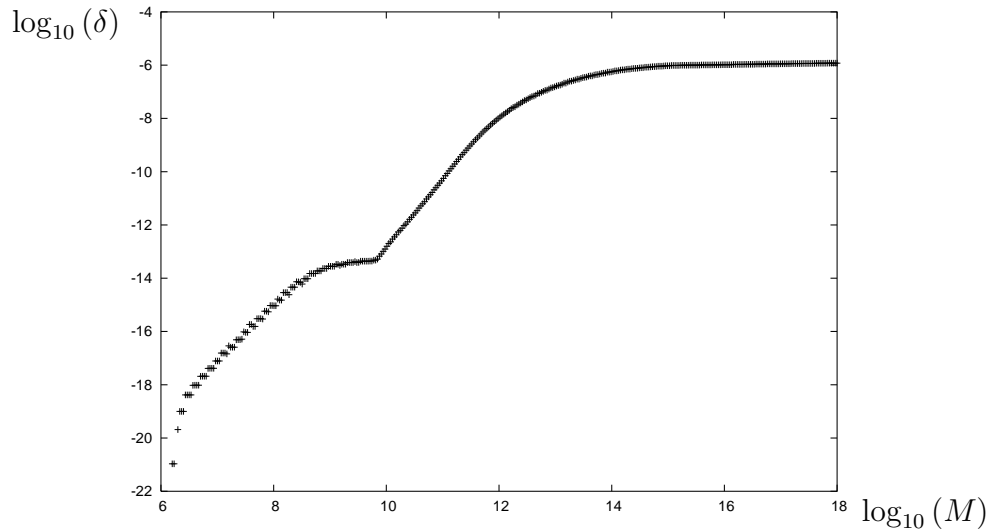}}}
    \put(125,0){$\log_{10}{(M)}$}
    \put(12,65){$\log_{10}{(\delta)}$}
  \end{picture}
  \caption{Log-log plot of $M$ against the largest value of $\delta$ found by
  a numerical computation.}\label{p:Mz2}
\end{figure}

As an illustration, we choose $M$ and the parameters $b_2$, $\gamma$,
$d$, $\nu$ based on a numerical experiment where we evaluate the
conditions of the theorem on a Cartesian grid of $360 \times 180
\times 120 \times 120 \times 20$ points in the cube $[10^6,10^{18}]
\times [10^{-6},1] \times [10^{-3}, 10] \times [10^{-3}, 10] \times
[10^{-3}, \frac{1}{9} \cdot 2^{-7/3}]$ in $(M,b_2,\gamma,d,\nu)$ space,
taking a logarithmic scale in the first four components and a linear
scale in the last.  For each $M$, we compute the largest value of the
coefficient $\delta$ over all grid points where all the conditions
hold in order to obtain a good choice of parameter values. Figure
\ref{p:Mz2} shows the largest $\delta$ one can obtain for a range of $M$.

\begin{cor}
  If $K$ has an invariant torus with frequency $K_0'(I_0)$, then
  the return map of $K$ on the surface of section $\xi=0$
  has an Aubry-Mather set with rotation number $\omega$ for
  any $\omega > K_0'(I_0)$.
\end{cor}

The corollary follows immediately from Aubry-Mather
theory \cite{AD83,Mat82,Mo86b} and the monotonicity of $K_0'$.
\change{Indeed, the Poincar\'e map $P_0$ corresponding to $K_0$ is given by
  \[ P_0:(\phi,I) \mapsto (\phi + K_0'(I), I)\,. \]
  From the monotonicity of $K_0'$, it follows that this is a twist map. The
  proof of the main result of Ref.~\cite{cny} shows that the Poincar\'e map $P$
  corresponding to $K$ also satisfies the twist condition. Alternatively, this
  also follows directly from the fact that the function $U$ is
  superquadratic.
}

\noteb{The paragraph that was here has been edited and changed into the following remark.}

\begin{rem}\label{t:measure}
  The translated curve theorem used (in \cite{cny}) to prove theorem
  \ref{t:main} shows persistence of invariant tori under relatively large
  perturbations (compared to other KAM or twist theorems), at the expense of
  excluding all but the ``most quasiperiodic'' frequencies $\omega$. If the
  Diophantine condition is relaxed to include a nonzero measure of
  frequencies, then one can similarly show that the region outside the innermost invariant torus contains a set of invariant tori of positive measure, where the measure converges exponentially fast to full measure as $I \to +\infty$. The remaining dynamics in this region is presumably
  chaotic, but due to the low dimension of the system, there can not be any
  Arnol'd diffusion. Thus, it is at worst bounded chaos. Among many papers on
  this subject, see e.g.\ \cite{Nei81,Poe82b,BNS03} for further details.

  By contrast, the region on the inside seems to exhibit chaotic dynamics,
  as can be seen in figures \ref{p:pp1} -- \ref{p:pp3}. In the near-autonomous
  setting (i.e., for small-amplitude $V$), the figures show large domains of
  integrable dynamics in the interior region.  A perturbation analysis for a
  similar system shows that, as the amplitude of $V$ goes to
  infinity, the sizes of these islands vanish, but their number increases
  reciprocally, so that an integrable set of measure $O(1)$ remains \cite{nst97}.
\end{rem}

In the next two sections, we consider potentials $V(x)$ with one and
two wave numbers.  The former case describes BECs in periodic lattices
\cite{bronski,promislow,pethick2,mapbecprl,mapbec}, and the latter case, which is now experimentally accessible \cite{quasibec}, describes BECs in periodic superlattices
\cite{eksioglu,louis2,super}.

\subsection{Example 1: BECs in Periodic Lattices}\label{lattice}

Optical lattice potentials are created experimentally as interference
patterns of counter-propagating laser beams \cite{olstandard}.  In the periodic case,
the external potential is typically taken to be sinusoidal,
\begin{equation}
        V(x) = -V_1\cos(\kappa x)\,,  \label{sing}
\end{equation}
where $\kappa = 2\pi/T$ is the lattice wave number.

We write equation (\ref{dynam35}) in the form
\begin{align}
        R' &= S\,, \notag \\    
        S' &= - \alpha_1 R + \alpha_3 R^3 + V_1 R \cos(\kappa x)\,,
 \notag \\
        x' &= 1 \,, \label{dyno}
\end{align}
where $\alpha_1 \propto \mu$, $\alpha_3 \propto g$, and $'$ denotes
differentiation with respect to $x$. The parameters $V_1\,, \alpha_1\,,
\alpha_3 \in \mathbb{R}$ and $\kappa > 0$ can all be adjusted
experimentally.  The phase variables are $(R,S) \in \mR^2$ and $x \in
\mR/(2\pi \kappa^{-1} \mZ)$.  The associated Poincar\'e map $P$, which
is defined to be the first return map on the section $x = 0$,
corresponds to the flow over $2\pi\kappa^{-1}$.

Equation (\ref{dyno}) has two reversible symmetries:
\[ (R,S,x) \mapsto (-R,S,-x) \,, \quad (R,S,x) \mapsto
(R,-S,-x). \] 
The map $P$ is thus reversible under reflection in both
coordinate axes,
        \[ R_j \circ P \circ R_j = P^{-1}\,, \quad j = 1,2\,, \]
where $R_1$ and $R_2$ are the reflections with respect to the axes.

Additionally, (\ref{dyno}) is invariant under two rescalings:
\begin{align}
  (R,S,x; \alpha_1,\alpha_3,\kappa,V_1) & \mapsto  \left(R,\lambda
  S,\lambda^{-1} x; \lambda^2 \alpha_1,\lambda^2 \alpha_3,\lambda
  \kappa,\lambda^2 V_1\right) \,, \label{e:s1} \\
  (R,S,x; \alpha_1,\alpha_3,\kappa,V_1) & \mapsto  \left(\mu R,\mu
  S,x; \alpha_1, \mu^{-2} \alpha_3,\kappa,V_1\right) \label{e:s2}\,.
\end{align}
The corresponding invariants for the Poincar\'e map are obtained by 
dropping the $x$ components.

One can rescale (\ref{dyno}) using these invariants to reduce to the
cases where $\alpha_1 = 0\,, \pm 1$ and $\alpha_3 = 0\,, \pm 1$.  We
look in detail at the case $\alpha_1 = -1$, $\alpha_3 = -1$,
corresponding to an attractive BEC with a negative chemical potential,
where the underlying integrable system (with $V_1 = 0$) is a ``figure
eight'' (consisting of a central saddle point and two exterior
centers, as shown in figure \ref{repulse1}c).

Because of the second rescaling (\ref{e:s2}), the parameter $\alpha_3$
simply measures the size of phase space.  The first rescaling
(\ref{e:s1}) shows that decreasing $\kappa$ increases the
nonintegrable perturbation by a square law.  Intuitively, a large
lattice wave number $\kappa$ implies that the Poincar\'e map
corresponds to short-``time'' integration in $x$ and is thus ``near''
the vector field.  With $\kappa$ rescaled to $1$, as was done in
stating and proving our main result, the perturbation $V(x)$ in
(\ref{dyno}) has unit period.

Figures \ref{p:pp1} -- \ref{p:pp3} show phase portraits of $P$ at
several amplitudes of the potential for
$(\alpha_1,\alpha_3,\kappa) = (-1,-1,1)$.  One obtains qualitatively
similar results for other values of the lattice depth $V_1$ if the
wave number $\kappa$ is rescaled, as indicated above. As remarked
previously, the phase space in these phase portraits is divided into
two clearly distinct regions: an outer one in which the dynamics
consists in large measure of invariant circles and Cantor-like
Aubry-Mather sets (that wind around the origin at large distance) and
an inner one in which the dynamics is mostly chaotic.  Our numerical
simulations, which show the same scaling that our theoretical results
indicate, suggest the presence of (parameter-dependent) integrable
dynamics of positive but small measure inside the ``chaotic sea.'' \cite{nst97}

\begin{rem}
  A similar combination of islands of invariant tori within a chaotic
  sea occurs in the example of a parametrically forced planar
  pendulum \cite{bhnsv}. The division of phase space into a mostly
  quasiperiodic and a mostly chaotic region is the typical behavior
  that one expects to observe in a large class of forced one {\it dof} Hamiltonian systems
  \cite{cny}.
\end{rem}


Applying theorem \ref{t:main} for $V(x)$ given by (\ref{sing}) implies
that the system (\ref{dyno}) has an invariant torus with frequency vector
$(\omega,1)$ provided $\omega$ satisfies the conditions of theorem
\ref{t:main}. The $R$-amplitude is roughly equal to $R_{\max}$, given by
\[ R_{\max} = 3 b_1 \frac{\hbar \kappa
  \omega}{\pi \sqrt{-mg}}\,, \] 
where $g < 0$ for the present case of attractive BECs [see equation \eq{I} and 
lemma \ref{t:K0}].


%
%


\subsection{Example 2: BECs in Periodic Superlattices}\label{superlattice}

Optical superlattices consist of small-scale lattices subjected to a
long-scale periodic modulation.  In recent experiments, BECs were
created in superlattice potentials with a length scale (wave number)
ratio of 1:3 \cite{quasibec}.  However, theoretical research concerning BECs in superlattices has only begun to gain prevalence \cite{eksioglu,louis2,super}. 

To consider the case of (symmetric) periodic superlattices, we examine
the potential
\begin{equation}
  V(x) = -\left[V_1\cos(\kappa_1 x) + V_2\cos(\kappa_2 x)\right]\,,  
\label{doub}
\end{equation}
where $\kappa_2 > \kappa_1$ without loss of generality and
$\kappa_2/\kappa_1 \in \mathbb{Q}$.  The minimal period is $T =
2\pi/\kappa$, where $\kappa := \mbox{gcd}(\kappa_1\,,\kappa_2)$.

Equation (\ref{dynam35}) is then written
\begin{align}
        \frac{dR}{dx} &= S\,, \notag \\    
        \frac{dS}{dx} &= - \alpha_1 R + \alpha_3 R^3 + V_1 R \cos(\kappa_1 x) 
+ V_2 R \cos(\kappa_2 x) \,, \notag \\
        \frac{dx}{dx} &= 1 \,, \label{dynsuper}
\end{align}
where all the parameters are again experimentally adjustable.


Applying theorem \ref{t:main} for $V(x)$ given by (\ref{doub}) with 
$\kappa_2/\kappa_1 \in \mQ$ (i.e., for periodic superlattices) implies that 
(\ref{dynsuper}) has an invariant torus with frequency vector $(\omega,1)$ 
provided $\omega$ satisfies the conditions of theorem
\ref{t:main}. As in the regular lattice case, the $R$-amplitude is roughly equal to
$R_{\max}$, where
\[ R_{\max} = 3 b_1 \frac{\hbar \kappa
  \omega}{\pi \sqrt{-mg}}\,, \] and we recall that $2\pi/\kappa$ is
the period of $V(x)$ and $g < 0$ for attractive BECs.

\begin{rem}
  If $\kappa_2/\kappa_1$ is not rational, then the potential $V$ is not periodic. If $\kappa_1$ and $\kappa_2$ satisfy a
  Diophantine condition, then one can prove the existence of invariant
  tori at large distance from the origin \cite{LZ95}.
  We conjecture that it is possible to quantify this existence result
  analogous to the periodic case \cite{cny}.
  If $\kappa_1$ and $\kappa_2$ are not Diophantine---for example, if
  $\kappa_2/\kappa_1$ is a Liouville number---then one expects
  unbounded solutions and no invariant tori \cite{hu98}.
\end{rem}


\section{Proof of the main result}\label{s:proof}

To prove Theorem \ref{t:main}, we first construct action-angle
coordinates explicitly. We then employ the KAM theorem of Chow, {\it et al.} \cite{cny} to complete the proof after a suitable transformation of the action variable.

\subsection{The action and $K_0$}

Define the action by
\[ I(h) = \int_{H_0 = h} S \, {\rm d} R = 4 \int_0^{R_{\max}} \sqrt{2 h
  - 2 z_4 R^4} \, {\rm d} R\,, \]
where $H_0(R,S) = \h S^2 + z_4 R^4$ and $R_{\max} = R_{\max}(h) > 0$ is
the solution of $H_0(R_{\max}, 0) = h$; that is, $R_{\max} =
z_4^{-1/4} h^{1/4}$.  With the substitution $u = R/R_{\max}$ and the
relation $h = z_4 R_{\max}^4$, the above integral reduces to
\begin{equation}
  I(h) = 4 \sqrt{2} b_1 z_4^{1/2} R_{\max}^3 = 4 \sqrt{2} b_1
  z_4^{-1/4} h^{3/4}\,. \label{e:I}
\end{equation}

\begin{lem}\label{t:K0}
  The unperturbed Hamiltonian in action-angle coordinates is given by
  \[ K_0(I) = 2^{-10/3} b_1^{-4/3} z_4^{1/3} I^{4/3}. \]
\end{lem}

The proof follows from the calculation above, noting that the
unperturbed Hamiltonian $K_0$ is the inverse of the function $I$
[because $H_0(R,S) = K_0(I)$ for $I = I(H_0(R,S))$]. The function $I$
is invertible because \[ \fd{I}{h} = 3 \sqrt{2} b_1 z_4^{-1/4} h^{-1/4} > 0. \]

\begin{rem}
  One could also define the action to be the area enclosed by the
  curve $\h S^2 + \bar{z}_2 R^2 + z_4 R^4 = h$, where $\bar{z}_2$ is
  the average of $z_2$, or even the area enclosed by $H = h$, where
  $\xi$ is considered as a parameter (that is, $\xi' = 0$ in the
  unperturbed system). In the latter case, the action and angle will
  depend on $\xi$.  Although these two approaches each leave a smaller term
  in the perturbation than our choice, and are therefore theoretically
  more pleasing, they lead to technical difficulties in the estimates
  we need to perform, as the expressions for the action and angle will
  involve elliptic integrals that depend on the phase variables.
\end{rem}


\subsection{The angle}

In the upper half plane ($S \geq 0$), we define the angle by
\[ \phi(h,R) = \left(\fd{I}{h}\right)^{-1} \int_0^R
\frac{1}{\bar{S}(h,w)} \, {\rm d} w ~ ~ (\mbox{mod } 1)\,, \]
where $\bar{S}(h,R) = \sqrt{2 h - 2 z_4 R^4}$ is the positive solution of
$H_0(R,\bar{S}(h,R)) = h$. Note that $\fd{I}{h} = \oint_{H_0 = h}
S^{-1} \, {\rm d} R$. A similar definition holds in the lower half
plane, where $\phi = \h - (\fd{I}{h})^{-1} \int_0^R \bar{S}^{-1} \,
{\rm d} w$ (mod 1). Henceforth, we restrict to the upper half plane
without loss of generality.

\begin{lem}The following formulas hold:\label{t:angle}
  \begin{eqnarray*}
    \phi(h,R) & = & \frac{1}{6 b_1} \int_0^{R/R_{\max}}
    (1-u^4)^{-1/2} \, {\rm d} u\,,\\
    \bar{S}(h,R) \fd{\phi}{h}(h,R) & = & -\frac{\sqrt{2}}{24 b_1} h^{-1/2}
    \frac{R}{R_{\max}}\,.
  \end{eqnarray*}
\end{lem}
Defining $r = r(h,R) = R/R_{\max}(h)$, it follows that
$\phi$ depends only on $r$:
\[ \phi(h,R) = \bar{\phi}(r(h,R)) \mbox{, where } \bar{\phi}(r)
= \frac{1}{6 b_1} \int_0^r (1-u^4)^{-1/2} \, {\rm d} u\,. \]
In particular, $\bar{\phi}(1) = 1/4$.

\begin{prf}
  Using the definition of $\phi$ and the substitution $w = R_{\max}
  u$, we obtain
  \begin{eqnarray*}
    \phi & = & \left(\fd{I}{h}\right)^{-1} \frac{1}{\sqrt{2}}
    z_4^{-1/2} R_{\max}^{-1} \int_0^r (1-u^4)^{-1/2} \,
    {\rm d} u\\
    & = & \frac{1}{6 b_1} z_4^{-1/4} h^{1/4} R_{\max}^{-1}
    \int_0^r (1-u^4)^{-1/2} \, {\rm d} u\,,
  \end{eqnarray*}
  which proves the first formula. Furthermore,
  \begin{eqnarray*}
    \bar{S}(h,R) \fd{\phi}{h} & = & -(2 h - 2 z_4 R^4)^{1/2}
    \frac{1}{6 b_1} (1-r^4)^{-1/2} R R_{\max}^{-2} \fd{R_{\max}}{h}\\
    & = & -\frac{\sqrt{2 h}}{24 b_1} R_{\max}^{-1} z_4^{-1/4} h^{-3/4}
    r\\
    & = & -\frac{\sqrt{2}}{24 b_1} h^{-1/2} r\,.
  \end{eqnarray*}
\end{prf}

\subsection{Localization and rescaling}

The system corresponding to $K$ can be computed directly and is
given by
\begin{eqnarray*}
  \phi' & = & K_0'(I) +f_1(\phi,I,\xi)\,,\\
  I' & = & f_2(\phi,I,\xi)\,,\\
  \xi' & = & 1\,,
\end{eqnarray*}
where $f_1(\phi,I,\xi) = - S \fd{\phi}{h} \fd{H_1}{R}$ and
$f_2(\phi,I,\xi) = -S \fd{I}{h} \fd{H_1}{R}$.

For a fixed $I_0$, we define a localization transformation $(\phi, I,
\xi) \mapsto (\phi, J, \xi)$ by $I = I_0 + \beta(I_0) J$, with
$\beta(I_0) = b_2 I_0 [\log(z_4 I_0)]^{-1}$, where $b_2 \in \mathbb{R}$ is positive and
will be determined later.  This takes the unperturbed torus $I = I_0$
to $J = 0$ and rescales the action variable, so that the components
$f_1$, $f_2$ of the perturbation are roughly the same size after
rescaling and some other conditions are met (see remark
\ref{t:rescaling} below). Although it is not a symplectic
transformation, it nonetheless maps Hamiltonian systems to Hamiltonian
systems.

Define $\omega = K_0'(I_0)$ and $m = \beta(I_0) K_0''(I_0)$, so that
\begin{eqnarray}
  \omega & = & \frac{1}{3} \cdot 2^{-4/3} b_1^{-4/3} z_4^{1/3}
  I_0^{1/3}\,, \label{e:omega}\\
  m & = & \frac{1}{9} \cdot 2^{-4/3} b_1^{-4/3} b_2
  \frac{z_4^{1/3} I_0^{1/3}}{\log{(z_4 I_0)}}\,. \label{e:m}
\end{eqnarray}

\begin{lem}\label{t:sys_rescaled}
  In $(\phi, J, \xi)$ coordinates, the system is written
  \[ \phi' = \omega + m J + g(J) + F_1(\phi,J,\xi) ~,~
  J' = F_2(\phi,J,\xi) ~,~ \xi' = 1, \]
  where, for some $J_*, J^* \in (0,J)$,
  \begin{eqnarray*}
    g(J) & = & \h \beta(I_0)^2 K_0'''(I_0+\beta(I_0) J_*)
    J^2\\
    & = & -\frac{1}{27} \cdot 2^{-4/3} b_1^{-4/3} z_4^{1/3} I_0^{1/3}
    \left(1+\frac{b_2 J_*}{\log{(z_4 I_0)}}\right)^{-5/3}
    \left(\frac{b_2 J}{\log{(z_4 I_0)}}\right)^2 \,, \\
    \fd{g}{J}(J) & = & \beta(I_0)^2 K_0'''(I_0+\beta(I_0) J^*) J\\
    & = & -\frac{1}{27} \cdot 2^{-1/3} b_1^{-4/3} b_2 \frac{z_4^{1/3}
      I_0^{1/3}}{\log{(z_4 I_0)}} \left( 1 + \frac{b_2 J^*}{\log{(z_4
    I_0)}} \right)^{-5/3} \frac{b_2 J}{\log{(z_4 I_0)}}\,, \\
    F_1(\phi, J, \xi) & = & f_1(\phi, I_0 + \beta(I_0) J, \xi)\\
    & = & \frac{1}{3} \cdot 2^{-2/3} b_1^{-2/3} z_2(\xi)
    r^2 z_4^{-1/3} I_0^{-1/3} \left(1+\frac{b_2 J}{\log{(z_4
        I_0)}}\right)^{-1/3}\,,\\
    F_2(\phi, J, \xi) & = & \beta(I_0)^{-1} f_2(\phi, I_0 +
    \beta(I_0) J, \xi)\\
    & = & -3 \cdot 2^{1/3} b_1^{1/3} z_2(\xi) r (1-r^4)^{1/2} b_2^{-1}
    \frac{\log{(z_4 I_0)}}{z_4^{1/3} I_0^{1/3}}
    \left(1+\frac{b_2 J}{\log{(z_4 I_0)}}\right)^{2/3}\,.
  \end{eqnarray*}
\end{lem}

The lemma shows that the nonintegrable parts $F_1$ and $F_2$ are order
$I_0^{-1/3}$ in leading term; this would have been order $1$ had the original function $U$ included a cubic term.

\begin{prf}
  Write $\Omega(J) = K_0'(I_0 + \beta(I_0) J)$, so that $\omega =
  \Omega(0)$ and $m = \Omega'(0)$.  By Taylor's theorem,
  \[ g(J) = \Omega(J) - \omega - m J = \h \Omega''(J_*) J^2 ~,~ g'(J)
  = \Omega'(J) - m = \Omega''(J^*) J\,. \]
  Furthermore, $\beta(I_0) J' = I'$. The expressions for $g$, $g'$,
  $F_1$, and $F_2$ follow by direct calculation from \eq{I} and
  lemma \ref{t:angle}.
\end{prf}

\subsection{Proof of Theorem \ref{t:main}}

The rescaled system is a perturbation of $\phi' = \omega + m J +
g(J)$, $J' = 0$. The unperturbed system has an invariant torus $J = 0$
with $\phi$-frequency $\omega$. Assume that $\omega$ is of constant
type with parameter $\gamma > 0$. We now study the persistence of this
torus under the perturbation $F = (F_1, F_2)$.

For $d > 0$, $\eta = 18 \gamma$, and $\rho = (3 m)^{-1} \eta$, we define
a (complex) neighborhood $\D_0$ of the unperturbed torus by
\[ \D_0 = \D_0(\eta,\rho, d) = \{ (\phi,J,\xi) \in \mC/\mZ \times
\mC \times \mR/\mZ: |\im(\phi)| \leq \eta, |J| < \rho \} + 2 d\,. \]
The KAM theorem of Chow, {\it et al.} \cite{cny} now shows that the perturbed
system has an invariant torus with frequency $\omega$ satisfying $|J|
\leq \frac{11}{19} d + \rho$ provided
\begin{eqnarray}
  \left |\left |\fd{g}{J}\right |\right |_{\D_0} & < & m/4\,, \label{e:twistcond}\\
  \omega_0 & \leq & \omega\,. \label{e:smallcond}
\end{eqnarray}
The first of these conditions is a twist condition, and the second
states that the perturbation is small enough, where $\omega_0$ is
defined as follows. Let $c = 5.04$,
and $||F||_{\D_0} = \max\{||F_1||_{\D_0}, ||F_2||_{\D_0}\}$. Let $L_W$
denote the Lambert $W$ function (i.e., the inverse of $W \mapsto W
e^W$).  Then,
\begin{eqnarray}
  \omega_0 & = & \alpha_0 \max\left\{1, \frac{1}{\log{2}} L_W(b
    \log{2})\right\}\,, \nonumber \\
  b & = & \frac{2+3 m}{c \gamma^2} \max\left\{12, \frac{7^8}{108
      (\eta-6 \gamma)^4}\right\} \max\{m d, 2 ||F||_{\D_0}\} \,,
  \label{e:b} \\
  \alpha_0 & = & \frac{3 + 12(\eta+2 d)}{2 d} \left( ||F||_{\D_0} + 2
    C \max\left\{ 1, \frac{2 ||F||_{\D_0}}{m d} \right\} \right) \,,
  \nonumber \\
  C & = & \frac{2}{3} \max\{ m (2 d + \rho) + ||g||_{\D_0} +
  ||F||_{\D_0}, 1\}\,. \nonumber
\end{eqnarray}

Lemmas \ref{t:twist} -- \ref{t:omega0} show that the twist and smallness
conditions follow from conditions \eq{cond1} -- \eq{cond_last} of
Theorem \ref{t:main}.

\begin{lem}\label{t:twist}
  Condition \eq{twistcond} follows from \eq{condL}.
\end{lem}

\begin{prf}
  For $|J| \leq \rho + 2 d$, it follows that
  \[ \lnorm \frac{b_2 J}{\log{(M)}} \rnorm_{\D_0} \leq
  b_2 \frac{\rho + 2 d}{\log{(M)}} = L \leq \frac{47}{200}\,. \]
  Hence, for $J^*$ as in lemma \ref{t:sys_rescaled},
  \[ \lnorm \left( 1 + \frac{b_2 J^*}{\log(M)} \right)^{-5/3}
  \frac{b_2 J}{\log(M)} \rnorm_{\D_0} \leq \left( 1 -
    \frac{47}{200} \right)^{-5/3} \frac{47}{200} < \frac{3}{8}\,. \]
  The desired result now follows from lemma \ref{t:sys_rescaled}.
\end{prf}

\begin{lem}\label{t:r}
  Let $r = R/R_{\max}$, as before.  Then,
  \begin{eqnarray*}
    ||r||_{\D_0} & \leq & 1+24 b_1 (\eta+ 2 d)\,, \\
    ||r (1-r^4)^{1/2}||_{\D_0} & \leq & [1+24 b_1 (\eta+ 2 d)]^3.
  \end{eqnarray*}
\end{lem}

\begin{prf}
  Without loss of generality, we restrict to $\phi$ satisfying $\re(\phi) \in
  [0,\frac{1}{4}]$ (mod 1). Observe that $|\im{\phi}| \leq \eta + 2 d$
  in $\D_0$. Thus, we can consider $\phi \in
  \bar{\D} = \{ \phi \in \mC/\mZ: \re(\phi) \in [0,\frac{1}{4}]
  (\mbox{mod } 1), |\im{\phi}| \leq \eta + 2 d \}$.

  \begin{figure}
    \begin{picture}(150,18)
      \put(75,9){\psarc(-1.7,0){0.8}{90}{270}}
      \put(75,9){\qline(-1.7,0.8)(1.7,0.8)}
      \put(75,9){\psarc(1.7,0){0.8}{270}{90}}
      \put(75,9){\qline(-1.7,-0.8)(1.7,-0.8)}
      \put(58,1){\circle*{1}}
      \put(58,17){\circle*{1}}
      \put(92,1){\circle*{1}}
      \put(92,17){\circle*{1}}
      \put(52.4, 8){1}
      \put(74.1, 3){2}
      \put(95.8, 8){3}
      \put(74.1, 13){4}
    \end{picture}
    \caption{The set $\E$, with labels for the four parts of its
      boundary (two semicircles and two line
      segments).}\label{p:stadium}
  \end{figure}
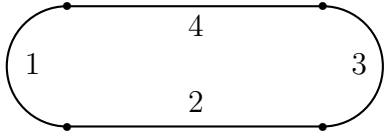

  Let $\E(\eps) = [0,1] + \eps \subset \mC$. It is obvious that $r$
  will be in $\E$ for large enough $\eps$. Let $\eps$ be the smallest
  number such that $r = r(\phi) \in \E = \E(\eps)$ for all $\phi \in
  \bar{\D}$.  Hence, there exists $\phi_0 \in \bar{\D}$ such that
  $r(\phi_0)$ is on the boundary of $\E$. Furthermore, $r(\re(\phi_0))
  \in [0,1]$.  Therefore, by the mean value theorem,
  \[ \eps \leq |r(\phi_0) - r(\re(\phi_0))| \leq
  \lnorm\fd{r}{\phi}\rnorm_{\bar{\D}} (\eta+2 d)\,. \]
  By definition, $\fd{\phi}{r}(r) = (6 b_1)^{-1} (1-r^4)^{-1/2}$.  Thus,
  \[ \fd{r}{\phi}(\phi) = 6 b_1 (1-r^4)^{1/2}\,. \]
  We will now show that $||1-r^4||_{\E} \leq (1+\eps)^4$.
  By the maximum modulus theorem, $||1-r^4||_{\E}$ is attained on the
  boundary of $\E$, which looks like a stadium and consists of four
  parts (see figure \ref{p:stadium}). Part 1 is parametrized by $r =
  \eps \exp(i \theta)$, where $\theta \in [\pi/2, 3 \pi/2]$. Hence,
  $|1-r^4| \leq 1+ \eps^4 \leq (1+\eps)^4$. Part 2 is parametrized by
  $r = u - \eps i$, with $u \in [0,1]$. Consequently, $|1-r^4| \leq 1-u^4 + 4
  u^3 \eps + 6 u^2 \eps^2 + 4 u \eps^3 + \eps^4 \leq 1+ 4 \eps + 6
  \eps^2 + 4 \eps^3 + \eps^4 \leq (1+\eps)^4$. Part 3 is parametrized
  by $r = 1 + \eps \exp(i \theta)$, where $\theta \in [-\pi/2,
  \pi/2]$. Hence, $|1-r^4| \leq 4 \eps + 6 \eps^2 + 4 \eps^3 + \eps^4
  \leq (1+\eps)^4$. Part 4 is similar to part 2.

  Thus, we have found that
  \[ \eps \leq 6 b_1 \left |\left |1-r^4\right |\right |^{1/2}_{\E} (\eta+2 d) \leq 6
  b_1 (\eta+2 d) (1+\eps)^2\,. \] 
  A straightforward calculation with
  this result then shows that $\eps \leq 24 b_1 (\eta + 2 d)$, which
  implies the first statement of the lemma. To prove the second
  statement, we use the fact that $||1-r^4||^{1/2}_{\D_0} \leq (1+\eps)^2 \leq
  [1 + 24 b_1 (\eta + 2 d)]^2$, as already shown.
\end{prf}

\begin{lem}\label{t:gF}
  Assume that \eq{condM}, \eq{condL}, \eq{condF}, and \eq{cond_last}
  hold. Then,
  \begin{eqnarray}
    ||g||_{\D_0} & \leq & \frac{2^{-4/3}}{27} b_1^{-4/3} (1-L)^{-5/3}
    L^2 M^{1/3}\,, \label{e:g} \\
    ||F||_{\D_0} & \leq & \nu d b_1^{-4/3} b_2
    \frac{M^{1/3}}{\log(M)}\,. \label{e:F}
  \end{eqnarray}
\end{lem}

\begin{prf}
  The assumption on $M$ implies that $\log{M} > 0$.  As in the proof of
  lemma \ref{t:twist}, $\lnorm \frac{b_2 J}{\log{(M)}} \rnorm_{\D_0}
  \leq L < 1$, and the result then follows from lemma
  \ref{t:sys_rescaled} and lemma \ref{t:r}. Condition \eq{condF}
  implies that $||F_1||_{\D_0} \leq ||F_2||_{\D_0}$, and
  \eq{cond_last} bounds $||z_2||_{\D_0}$.
\end{prf}

\begin{lem}\label{t:omega0}
  Assume that $\nu \leq \frac{1}{9} \cdot 2^{-7/3}$ and $\gamma =
  \eta/18 \leq \frac{49}{72}$, as in Theorem \ref{t:main}.
  If conditions \eq{condM} -- \eq{cond_last} hold, then
  \[ \omega_0 \leq \frac{A}{\log(2)} \left(\frac{1}{3} +
    \frac{\log(B)}{\log(M)}\right) M^{1/3}\,. \]
\end{lem}
Theorem \ref{t:main} now follows immediately from this lemma and the
expression for $\omega = K'(I_0)$.

\begin{prf}
  By the previous lemma,
  \[ C \leq \frac{2}{3} \max\left\{\frac{\eta}{3} + \left(
      \frac{2^{-1/3}}{9} + \nu\right) d b_1^{-4/3} b_2
    \frac{M^{1/3}}{\log(M)} + \frac{2^{-4/3}}{27} b_1^{-4/3}
    (1-L)^{-5/3} L^2 M^{1/3} , 1\right\}\,. \]
  A straightforward calculation shows that $\nu \leq \frac{1}{9} \cdot
  2^{-7/3}$ implies $2 ||F||_{\D_0} \leq m d$. Hence,
  \begin{eqnarray*}
    \alpha_0 & = & \frac{3 + 12 (\eta+2 d)}{2 d} (||F||_{\D_0} + 2 C)\\
    & \leq & \frac{1 + 4 (\eta+2 d)}{d}
    \max\left\{ \frac{2 \eta}{3} M^{-1/3} \log(M) +
      \left( \frac{2^{2/3}}{9} + \frac{7 \nu}{2}\right) d b_1^{-4/3}
      b_2 +\right.\\
    & & \left. \frac{2^{-1/3}}{27} b_1^{-4/3} (1-L)^{-5/3} L^2 \log(M),
      2 M^{-1/3} log(M) + \frac{3 \nu}{2} d b_1^{-4/3} b_2 \right\}
    \frac{M^{1/3}}{\log(M)}\\
    & = & A \frac{M^{1/3}}{\log(M)}\,.
  \end{eqnarray*}
  From $\gamma = \eta/18 \leq 49/72$, it follows
  that $12 \leq \frac{1}{108} \cdot 7^8 (\eta-6\gamma)^{-4}$. We observe
  that $\gamma$ does not appear in any other quantity besides $b$. The
  choice $\gamma = \eta/18$ minimizes $b$. We obtain
  \begin{eqnarray*}
    b & = & \frac{7^8 \cdot 18^2 \cdot 3^4 (2 + 3 m)}{108 \cdot 2^4 c
      \eta^6} m d\\
    & = & 2^{-4} \cdot 3^5 \cdot 7^8 c^{-1} \eta^{-6} (2+3 m) m d\\
    & = & 2^{-13/3} \cdot 3^3 \cdot 7^8 c^{-1} d b_1^{-4/3} b_2 \eta^{-6}
    \left(\frac{1}{\log(M)} + \frac{2^{-7/3}}{3} b_1^{-4/3} b_2
      \frac{M^{1/3}}{\log^2(M)} \right) M^{1/3}\\
    & = & B M^{1/3}\,.
  \end{eqnarray*}
  From $b \geq 2$, we see that $L_W(b \log(2)) \leq \log{(b)}$. The
  result follows.
\end{prf}

\begin{rem}\label{t:rescaling}
  The transformation $I \mapsto J$ is chosen so that $\omega$ and
  $\omega_0$ have the same growth rate in $I_0$ and $z_4$ (i.e., in
  the leading term).  Independent of the rescaling, $\omega$ grows as
  $I_0^{1/3}$. On the other hand, $\omega_0 \leq \alpha_0
  \log{(b)}/\log{(2)}$ (if $b \geq 2$). Because $b$ grows as some
  power of $I_0$, it follows that $\alpha_0$ must grow as $I_0^{1/3} /
  \log{(I_0)}$.  Analyzing the dependence on $z_4$ then leads to the
  chosen rescaling.
\end{rem}

\section{Summary and Conclusions} \label{conc}

In this paper, we discussed quasiperiodic dynamics in Bose-Einstein condensates (BECs) in periodic lattices and superlattices.  The mean-field dynamics of a BEC is governed by the Gross-Pitaevskii (GP) equation (\ref{GPE}), which consists of a cubic nonlinear Schr\"odinger equation plus an external potential that takes into account the \textquotedblleft trap" where the condensate resides.  In this mean-field description, a given particle in the BEC is affected by the other particles only through average effects, and the non-local term in the original many-body Hamiltonian leads to the nonlinear term in the GP equation.  One obtains a cubic nonlinearity if considering only two-body interactions.  The GP equation, which is derived as a zero-temperature theory, provides a good description for BEC dynamics below the critical transition temperature at which the condensate forms \cite{stringari,lieb}.  In \textquotedblleft cigar-shaped" BECs, two dimensions are tightly confined, so one may further reduce (\ref{GPE}) to (\ref{nls3}), which has one spatial dimension \cite{stringari}.  The amplitude dynamics of coherent structures of (\ref{nls3}) with trivial phase (which describe standing waves) are governed by a forced Duffing equation given by (\ref{dynam35}), or equivalently by \eq{sys}. The first equation arises directly from the physical setting, whereas \eq{sys} is the more convenient decription for mathematical analysis. We briefly describe the dynamics of this system.

In the absence of the forcing (i.e., when the external potential $V \equiv 0$), the system reduces to the
autonomous Duffing oscillator and is integrable. Its dynamics depends on the
sign of the chemical potential $\mu$ (which indicates how many particles are trapped in the BEC) and is illustrated in figure \ref{repulse1}.  Here, we considered the two situations with negative chemical potential (figure \ref{repulse1}b,c).  In both cases, there is a family of invariant tori
winding around the trivial periodic orbit $R \equiv 0$. When viewed in the
proper coordinates, the dynamics on each torus is a flow with constant velocity
vector.  One frequency is equal to 1, and the other varies from one torus
to the next, going to infinity monotonically as the distance of the torus to
$R \equiv 0$ goes to infinity. The ``proper coordinates'' are the action-angle
coordinates introduced in section \ref{s:proof}.  These tori correspond to
quasiperiodic oscillations of $R(x)$, which describes the amplitude dynamics of the BEC wave function. We note that in this unforced setting, the quasiperiodicity has no physical meaning, as the period of
the optical lattice potential does not yet play a role.

In the forced setting, $V \neq 0$.  Figures \ref{p:pp1}
-- \ref{p:pp3} illustrate the case of a periodic lattice for a range of
amplitudes. We remark that from a mathematical point of view, there is no
extra difficulty in dealing with superlattices or even lattices with
more than two Fourier modes as compared to lattices with just a single mode,  provided the periods are commensurate (that is, provided they have a
common multiple). Our main result, Theorem \ref{t:main}, uses KAM theory to
demonstrate that for any size of the periodic forcing, there exist invariant
tori winding around the trivial periodic orbit. Let $R_{\max}$ denote the
maximal $R$ coordinate on an invariant torus.  That is, $R_{\max}$ is the
amplitude of the corresponding oscillations. In terms of the mathematical
setting \eq{sys}, invariant tori exist for a dense set of ``admissible''
$R_{\max}$ in the interval $(c z_4^{-1/2} ||z_2||^{1/2}, +\infty)$, for some
positive constant $c > 0$. In terms of the physical parameters, the lowest
admissible amplitude is roughly equal to $3 b_1 \hbar \kappa \omega/(\pi
\sqrt{-m g})$, where $\omega$ is its frequency, $\kappa$ is the lattice
wave number, $b_1 \approx 0.874019$ is a constant, $m$ is the mass of the atomic species in the condensate, and $g < 0$ is the scaled value of the two-body scattering length.  The parameter $g$ varies from one condensate species to another and can be changed by exploiting Feshbach resonances \cite{fesh}.

\noteb{The following paragraph has been edited.}

The set of admissible $R_{\max}$ is characterized by a Diophantine condition
on the frequencies of the torus. A slightly different KAM theorem than the one used in this work shows that, for each forcing, the relative measure of invariant
tori converges to full measure as $R_{\max} \to +\infty$.  See, for example,
\cite{Nei81,Poe82b} for a general statement or \cite{BNS03} for discussion of a
system similar to ours.  Thus, we find a large measure of quasiperiodic
dynamics, and in the complement of the union of invariant tori, there exist
Aubry-Mather sets. Furthermore, by the Poincar\'e-Birkhoff Theorem
\cite{Birk13,Birk25}, in between any two invariant tori there are periodic
orbits of saddle and center type for all intermediate resonant frequencies
(technically, these can be classified as Aubry-Mather sets as well).
As a result, there is a large measure of quasiperiodic invariant tori interlaced with
Aubry-Mather sets and resonant layers.
We further note that the homoclinic figure eight that exists in the unforced
setting for $\mu < 0$ will generically break under the forcing, generating a
homoclinic tangle, as can be checked for this specific case by computing
Melnikov integrals.  We refer to \cite{gucken,Wig88} for such computations on
similar systems.  A possible direction for future study would be to establish
conditions for the nonexistence of invariant tori in the system under
consideration, which would complement the existence conditions presented in
this paper. Such a \emph{converse KAM theorem} was obtained in \cite{Mat84}
for area-preserving twist maps and in \cite{MacKay89} in a Lagrangian/Hamiltonian setting.

From a physical perspective, we recall that the GP equation is derived from a many-body quantum problem under the assumption that the Bose gas is dilute \cite{stringari,lieb}.  In particular, this implies that the mean interparticle distance should be much larger than the scattering length.  Therefore, the large $R$ results obtained here would necessitate the condensate to have expanded sufficiently to ensure that its density is low (so that the reduction from the many-body problem to the three-dimensional GP description remains valid).  Accordingly, for sufficiently large $R$, one would eventually have to include corrections to the GP equation arising beyond the mean-field description to describe the physics correctly.  Incorporating beyond-mean-field dynamics in the study of BECs is a difficult problem (see, e.g., the discussion in Ref.~\cite{anamaria}), and it is not agreed precisely when such corrections become relevant or what form they should take to simultaneously ensure tractability and correctly capture the physics.

\section*{Acknowledgements}

We are grateful to Peter Engels, Panos Kevrekidis, Boris Malomed, Alexandru Nicolin, and Li You for useful discussions concerning this research.  We also thank the editor and an anonymous referee whose constructive comments lead to significant improvement of this manuscript.  MAP was supported in part by a VIGRE grant awarded to the School of Mathematics at Georgia Tech and in part by the Gordon and Betty Moore Foundation through Caltech's Center for the Physics of Information.  MvN was partially supported by the Center for Dynamical Systems and Nonlinear Studies at Georgia Tech and partially by EPSRC grant GR/S97965/01. YY was partially supported by NSF grant DMS0204119.


\end{document}